\documentclass[11pt,a4paper]{article}
\usepackage[numbers,sort&compress]{natbib}
\usepackage[colorlinks,
            linkcolor=blue,
            anchorcolor=green,
            citecolor=magenta
            ]{hyperref}
\usepackage{mathrsfs,amssymb,amsfonts}

\usepackage{graphicx}
\usepackage{amsmath,amssymb,latexsym,color}
\usepackage[mathscr]{eucal}
\usepackage[numbers]{natbib}
\usepackage{bm}
\usepackage{subfigure}
\usepackage[top=1in, bottom=1in, left=1.25in, right=1.25in]{geometry}
\usepackage{indentfirst}
\usepackage{array}
\usepackage{graphicx}
\usepackage{setspace}
\usepackage{epstopdf}

\newtheorem{theorem}{Theorem}[section]

\newtheorem{lemma}[theorem]{Lemma}
\newtheorem{remark}[theorem]{Remark}

\newtheorem{corollary}[theorem]{Corollary}
\newtheorem{example}{Example}

\newtheorem{question}[theorem]{Question}

\newcommand{\qed}{\hfill\Box\medskip}
\usepackage{CJK}
\usepackage{geometry}
\usepackage[authormarkup=none]{changes}
\definechangesauthor[name={Reviewer 1}, color=orange]{Reviewer 1}
\definechangesauthor[name={Reviewer 2}, color=purple]{Reviewer 2}  %修改
\begin{document}

\setlength{\baselineskip}{14pt}
\renewcommand{\abovewithdelims}[2]{
\genfrac{[}{]}{0pt}{}{#1}{#2}}
%%%%%%%%%%%%%%%%%%%%%%%%%%%%%%%%%%%%%%%%%%%%%%%%%%%%%%%%%%%%%%%%%%%%%%%%%%%%%%%%%%%%%%%%
%%%%%%%%%%%%%%%%%%%%%%%%%%%%%%%%%%%%%%%%%%%%%%%%%%%%%%%%%%%%%%%%%%%%%%%%%%%%%%%%%%%%%%%%

\title{\bf {A note on the second-largest number of dissociation sets in connected graphs} \footnote{This research is supported by the National Natural Science Foundation of China(12301453, 12271524, 12331013), the Project of Education Department of Hunan Province(23B0125), the Project of Scientific Research Fund of Hunan Provincial Science and Technology Department (2018WK4006), the Natural Science Foundation of  Jiangxi(20224ACB201002) and the Natural Science Foundation of Hunan (2022JJ30674).}}

\author{ Pingshan Li, \quad Ke Yang\footnote{Corresponding author. \newline {\em E-mail address:} lips@xtu.edu.cn (P. Li), jinweipei82@163.com(W. Jin), 202421511255@smail.xtu.edu.cn(K. Yang)}, \quad Wei Jin\\
{\textsuperscript{a}\footnotesize  School of Mathematics and Computational Science, Xiangtan University, Xiangtan, Hunan 411105, PR China}\\
{\textsuperscript{b}\footnotesize Key Laboratory of Intelligent Computing$\And$Information Processing of Education}\\
{\textsuperscript{c}\footnotesize Key Laboratory for Computation and Simulation in Science and Engineering}\\
{\textsuperscript{d}\footnotesize National Center for Applied Mathematics in Hunan}\\
}

\date{}
 \maketitle
\begin{abstract}
\setlength{\parindent}{2em}

A subset of vertices is called a dissociation set if it induces a subgraph with vertex degree at most one.
Recently, Yuan et al. established the upper bound of the maximum number of dissociation sets among all connected graphs of order $n$ and characterized the corresponding extremal graphs. They also proposed a question regarding the second-largest number of dissociation sets among all connected graphs of order $n$ and the corresponding extremal graphs.
In this paper, we give a positive answer to this question.
\medskip


\medskip
%\noindent {\em 2010 MSC:} 05C25; 05C15.
\end{abstract}

\section{Introduction}\label{section1}
\setlength{\parindent}{2em}

We consider finite, simple, and undirected graphs and follow the terminology in  \cite{bondy2008}.
The study of graph substructures and their enumeration has been a theme in extremal graph theory, dating back to a question posed by Erd\H{o}s and Moser in 1960:
what is the maximum possible number of maximal independent sets in general graphs of order $n$, and what are the graphs that achieve this maximum value?
This was solved by Moon and Moser \cite{moon} in 1965. Since then, the maximum number of maximal or maximum independent sets on various graph families such as trees, unicyclic graphs, connected graphs and the graphs with at most $r$ cycles have been determined \added[id=Reviewer 2]{; see \cite{wilf1986, griggs1988, KohGohDong2008, SaganVatter2006, Wloch2008, zito1991} for details}.

A subset $S$ of vertices in a graph $G$ is called a dissociation set if it induces a subgraph with vertex degree at most one. This concept, introduced by Yannakakis \cite{yannakakis1981}, not only generalizes the notions of independent sets and induced matchings but also forms a dual relationship with the 3-path vertex cover problem \cite{bresar2011}. In recent years, scholars have begun to pay attention to determining the maximum (or minimum) number of maximal or maximum dissociation sets on various graph families and characterize the \added[id=Reviewer 1]{corresponding} extremal graphs such as trees \cite{tu2021,TU2024,Wang2024}, forests \cite{SunLi2023}, unicyclic graphs \cite{Zhang2024(1)}, general and triangle-free graphs \cite{Tu2022(1)}.

Recently, Yuan et al. \cite{yuan2024} focused on the problem of the maximum number of all dissociation sets and the characterization of extremal graphs. They established  \added[id=Reviewer 1]{the maximum number} of dissociation sets among all connected graphs of order $n$ and characterized the corresponding extremal graphs. They also proposed a question  \added[id=Reviewer 1]{ of determining} the second-largest number of dissociation sets and the corresponding extremal graphs as follows.

\begin{figure}[htbp!]
  \centering
  \includegraphics[width=0.4\textwidth]{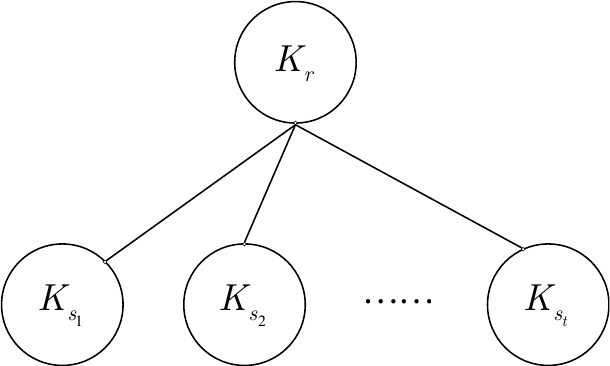}\\
  \caption{ $K_{r} \ast (K_{s_{1}}\cup K_{s_{2}}\cup \cdots \cup K_{s_{t}})$.}
  \label{figure Kr}
\end{figure}

\begin{question}\label{question1.1}\rm
\emph{Is $h(n)$ the second-largest number of dissociation sets among all connected graphs of order $n$ when $n\ge 10$? Do the extremal graphs belong to
\[\begin{cases}
  K_3*\frac{n-3}{2} K_2\ or\ K_1*(2K_1\cup \frac{n-3}{2}K_2)& \text{ if $n$ is odd; }  \\
  K_3*(\frac{n-4}{2} K_2\cup K_1)\ or\ K_1*(3K_1\cup \frac{n-4}{2}K_2)& \text{ if $n$ is even.}
\end{cases}\]} 
\end{question}
\added[id=Reviewer 1]{Here, $h(n)$ is defined as follows\[h(n) :=
\begin{cases}
2^{n-1} + (n+9) \cdot 2^{\frac{n-7}{2} }  & \text{if } n \text{ is odd}; \\
2^{n-1} + (n+12) \cdot 2^{\frac{n-8}{2} }  & \text{if } n \text{ is even and } n \neq 6; \\
42 & \text{if } n = 6.
\end{cases}\]}
The operation $*$ is a method of transforming several disjoint complete subgraphs into a connected graph, defined as follows. Let $G=K_{s_{1}}\cup K_{s_{2}}\cup \cdots \cup K_{s_{t}}$ be a disconnected graph whose all components are complete subgraphs, and $K_r$ a complete graph disjoint from $G$.  $K_r * G$ is the graph obtained by selecting a vertex $u \in V(K_r)$ and connecting $u$ to exactly one vertex from each component of $G$ (see Figure \ref{figure Kr}). 

In this paper, we give a positive answer to this question. The remainder of the paper is organized as follows: Section \ref{section2} introduces key definitions and preliminary lemmas; Section \ref{section3} presents the main results, including \added[id=Reviewer 1]{the second-largest number of dissociation sets in connected graphs of order $n$} and the corresponding extremal graphs.

\section{Preliminaries}
\label{section2}

The set of all dissociation sets of $G$ is denoted by $\mathscr{D}(G)$, and the number of dissociation sets consisting of $k$ vertices is denoted by $d\left ( G,k \right )$. For convenience, we define $d\left ( G,0 \right ) = 1$ for all graphs. The total number of dissociation sets $d(G)$ of $G$ is defined as \[d(G)=\sum_{k=0}^{\left | V(G )\right | } d(G,k).\]
Let $G_1, G_2, \dots, G_t $ be all components of $G$. We can easily see that $d(G) = \prod_{i=1}^{i=t} d(G_i)$. Two different vertices $u,v\in V(G)$ are called \emph{true twins} if $N\left [ u \right ] =N\left [ v \right ] $ which indicates that $uv\in E(G)$ and are called \emph{false twins} if $N(u)=N(v)$. It is easy to see that if \added[id=Reviewer 2]{$u$ and $v$} are true twins in graph $G$, then \added[id=Reviewer 2]{$u$ and $v$} are false twins in the subgraph $G-uv$.

A vertex with degree one is called a \emph{pendant vertex} and the neighbor of a pendant vertex with degree at least two is called a \emph{quasi-pendant vertex}. The dimension of cycle space \added[id=Reviewer 2]{of a connected graph $G$} is denoted by $c(G)=|E(G)|-|V(G)|+ 1$ where $|E(G)|$ and $|V(G)|$ represent the number of edges and vertices in $G$, respectively.
Consider the cycle $C_n$ as $v_1v_2\cdots v_nv_1$. Denote by $C_n(k_1,\cdots, k_n)$ the unicyclic graph obtained from $C_n$ by identifying $v_i\in V(C_n)$ with a vertex of pendant tree $T_{k_i}$ for $1\le i\le n$, where $k_1+\cdots +k_n$ equals to the order of \added[id=Reviewer 2]{$C_n(k_1,\cdots, k_n)$}.

\begin {lemma}\label{lemma2.1} \rm( \cite{yuan2024}).
\emph{Let $G$ be a graph of order $n$, then $d(G)\le 2^{n} $, the equality holds if and only if $G\cong sK_{1} \cup tK_{2} $ where $s+2t=n$}.
\end {lemma}

%In the content below, we set $g(n): = 2^{n}$.

\begin {lemma}\label{lemma2.2} \rm( \cite{yuan2024}).
\emph{Let $G$ be a graph with $v\in V(G)$. Then
\[d(G)=d(G-v)+d(G-N[v])+\sum_{u\in N(v)}^{} d(G-N[u]\cup N[v])\textcolor{purple}{.}\]}
\end {lemma}

\begin {lemma}\label{lemma2.3} \rm( \cite{yuan2024}).
\emph{Let $G$ be a graph with \added[id=Reviewer 2]{$uv\in E(G)$}. Then $d(G)\le d(G-uv)$, the equality holds if and only if $N_{G} [u]=N_{G} [v]$.}
\end {lemma}

\begin {theorem}\label{theorem2.4} \rm( \cite{yuan2024}).
\emph{Let $G$ be a connected graph of order $n$. Then $d(G)\le f(n)$, the equality holds if and only if $G\cong F_{n}$, where\[f(n) :=
\begin{cases}
2^{n-1} + (n+3) \cdot 2^{\frac{n-5}{2} }  & \text{if } n \text{ is odd}; \\
2^{n-1} + (n+6) \cdot 2^{\frac{n-6}{2} }  & \text{if } n \text{ is even,}
\end{cases}\]and \[F_n: = \begin{cases}
  K_{1} \ast \frac{n-1}{2} K_{2}& \text{ if $n$ is odd; } \\
  P_{6} \ or \ K_{2}\ast 2K_{2}   & \text{ if $n$ = 6; }  \\
  K_{2} \ast \frac{n-2}{2} K_{2}& \text{ if $n$ is even and $n$ $\ne$ 6.}
\end{cases}\]}
\end {theorem}

 \begin {theorem}\label{theorem2.5} \rm( \cite{yuan2024}).
\emph{Let $G$ be a unicyclic graph of order $n\ge 3$. Then $d(G)\le h(n)$, the equality holds if and only if $G\cong U_{n} $, where\[h(n) :=
\begin{cases}
2^{n-1} + (n+9) \cdot 2^{\frac{n-7}{2} }  & \text{if } n \text{ is odd}; \\
2^{n-1} + (n+12) \cdot 2^{\frac{n-8}{2} }  & \text{if } n \text{ is even and } n \neq 6; \\
42 & \text{if } n = 6 \textcolor{purple}{,}
\end{cases}\]and \[U_n := \begin{cases}
  K_{3}\ast (\left \lfloor \frac{n-3}{2}  \right \rfloor K_{2}\cup (n+1 \bmod 2 )K_{1}  ) & \text{ if } n\ne 6; \\
  K_{1}\ast (K_{3} \cup K_{2} ) & \text{ if } n=6.
\end{cases}\]}
\end {theorem}

\begin{lemma}\label{lemma2.6} Let $G$ be a graph \added[id=Reviewer 2]{of order $n\ge 4$} with quasi-pendant vertex $u_q$. Suppose that $v_1, v_2, \cdots, v_s$ are pendant vertices that are  adjacent to $u_q$ where $s\ge 2$. Then $d(G)<d(G+v_{s-1}v_s-u_qv_s)$.
\end{lemma}

\noindent{\bf Proof.} Note that \added[id=Reviewer 2]{$v_{s-1}$ and $v_{s}$} are true twins in $G+v_{s-1}v_s$.  By Lemma \ref{lemma2.3}, $d(G)=d(G+v_{s-1}v_{s})$. Note that \added[id=Reviewer 2]{$u_{q}$ and $v_{s}$} are not true twins in the graph $G+v_{s-1}v_{s}$. By Lemma \ref{lemma2.3}, $d(G+v_{s-1}v_{s})<d(G+v_{s-1}v_s-u_qv_s)$. Then $d(G)<d(G+v_{s-1}v_s-u_qv_s)$.
\hfill$\Box$

\noindent{\bf Remark.} The proof of this lemma comes from \cite{yuan2024}, but its conclusion is stronger than the lemma above. For the sake of the completeness of the paper and to better serve this paper, we have rewritten the conclusion above.

\begin {corollary}\label{corollary2.7} \rm
\emph{Let $T$ be a tree of order $n$ with the second-largest number of dissociation sets and $u_{q}$ a quasi-pendant vertex of $T$. If $u_q$ is connected to $s$ pendant vertices, then $s\le 3$.}
\end {corollary}

\added[id=Reviewer 1]{\noindent{\bf Proof.} In contrast, suppose that $s\ge 4$. Choose four pendant vertices $v_{1} ,v_{2},v_{3} ,v_{4}$ that are adjacent to $u_q$. By Lemma \ref{lemma2.6}, $d(T)<d(T+v_{1}v_{2}-u_qv_2):=d(T')$. Similarly, $d({T}')<d({T}'+v_{3}v_{4}-u_qv_4):=d(T'')$. Consequently, $d(T)<d({T}')<d(T'')$ and $|V(T)|=|V(T')|=|V(T'')|$. This contradicts the fact that $T$ is the tree of order $n$ with the second-largest number of dissociation sets. So, $s\le 3$.}
%\noindent{\bf Proof.} In contrast, suppose that $s\ge 4$. Choose four pendant vertices $v_{1} ,v_{2},v_{3} ,v_{4}$ that are adjacent to $u_q$. Note that $v_{1} ,v_{2}$ are true twins in $T+v_1v_2$.  By Lemma \ref{lemma2.3},  $d(T)=d(T+v_{1}v_{2})$.  Let ${T}'=T+v_{1}v_{2}-u_{q}v_{2}$ (see Figure \ref{figure s<3}). Note that $u_{q},v_{2}$ are not true twins in the graph $T+v_{1} v_{2}$. By Lemma \ref{lemma2.3}, $d(T+v_{1}v_{2})<d(T+v_{1}v_{2}-u_qv_2)=d(T')$. Similarly, we can find that $d({T}')=d({T}'+v_{3} v_{4})<d({T}'+v_{3}v_{4}-u_qv_4)$. Let $T''={T}'+v_{3}v_{4}-u_qv_4$.
%Consequently, $d(T)<d({T}')<d(T'')$ and $|V(T)|=|V(T')|=|V(T'')|$. This contradicts the fact that $T$ is the tree of order $n$ with the second-largest number of dissociation sets. So, $s\le 3$.
\hfill$\Box$

\begin{figure}[htbp!]
  \centering
  \includegraphics[width=0.9\textwidth]{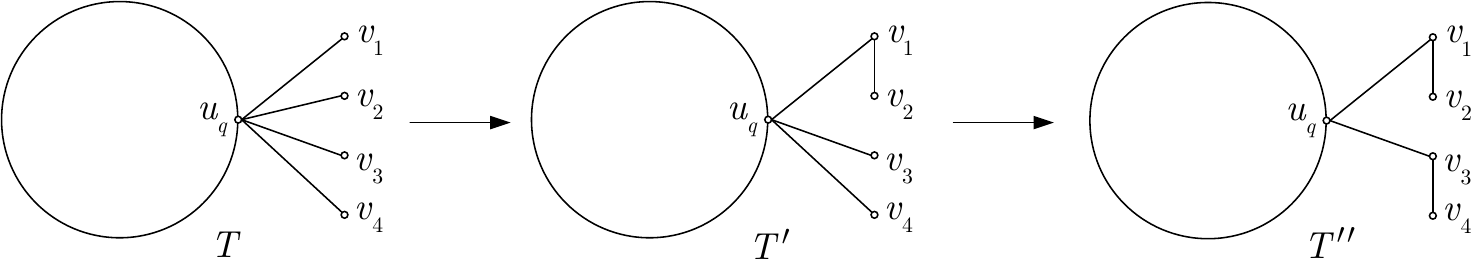}\\
  \caption{ $T$, ${T}'$, ${T}''$.}
  \label{figure s<3}
\end{figure}

\section{The second-largest number of dissociation sets in connected graph}
\label{section3}

In this section, we will give a positive answer to Question \ref{question1.1}. In fact, we can replace the prerequisite $n\ge 10$ with $n\ge 9$. 
First, we determine the second-largest number of dissociation sets among all trees of order $n$ and characterize the structure of the corresponding extremal trees.

\begin{lemma}\label{lemma3.1}
Let $T$ be a tree of order $n(n=9\ or\ 10)$ where $T\notin\{K_1\ast 4K_2,  K_2\ast 4K_2\} $. Then $d(T)\le h(|V(T)|)$, the equality holds if and only if $T\in\{K_1\ast (2K_1\cup 3K_2),  K_1\ast (3K_1\cup 3K_2)\}$ where $h(9)=292, h(10)=556$.

\end{lemma}

\noindent{\bf Remark.} The proof of this lemma is obtained by analyzing all trees of order $9$ or $10$. The conclusion is essential, but verification is unskilled and very cumbersome. So we provide the detailed proof in the appendix.

\begin{theorem}\label{theorem3.1}
Let $T$ be a tree of order $n (n\ge 9)$. If $T\ncong F_{n} $, then $d(T)\le h(n)$, the equality holds if and only if $T\cong T_{n} $, where \[T_{n} : =\begin{cases}
  K_{1}\ast(2K_{1}\cup\frac{n-3}{2}K_{2}),& \text{ if $n$ is \textcolor{purple}{odd;}}\\
  K_{1}\ast(3K_{1}\cup\frac{n-4}{2}K_{2}),& \text{ if $n$ is \textcolor{purple}{even,}}
\end{cases}\]$F_{n} $ and $h(n)$ are defined in Theorem \ref{theorem2.4} and Theorem \ref{theorem2.5}, respectively.
\end{theorem}

\noindent{\bf Proof.} We shall prove by induction on $n$. By Lemma \ref{lemma3.1}, the theorem holds for $n=9$ and $n=10$. Assume that the theorem is true when the order of the tree is less than $n$ where $n\ge11$. Let $T$ be a tree of order $n$ such that $T\ncong F_n$. By Theorem \ref{theorem2.4}, $d(T)<f(n)$. Let $T$ be a tree of order $n$ with the second-largest number of dissociation sets and choose a longest path $P=v_{0}v_{1}v_{2}\cdots v_{l} $ of $T$ starting from $v_{0}$.  By Corollary \ref{corollary2.7}, we can obtain ${\rm deg}_T(v_1)\le 4$. To be more specific, we have the following claim.

{\bf Claim 1. }  $ {\rm deg}_T(v_{1})=2$.

\begin{figure}[htbp!]
  \centering
  \includegraphics[width=0.6\textwidth]{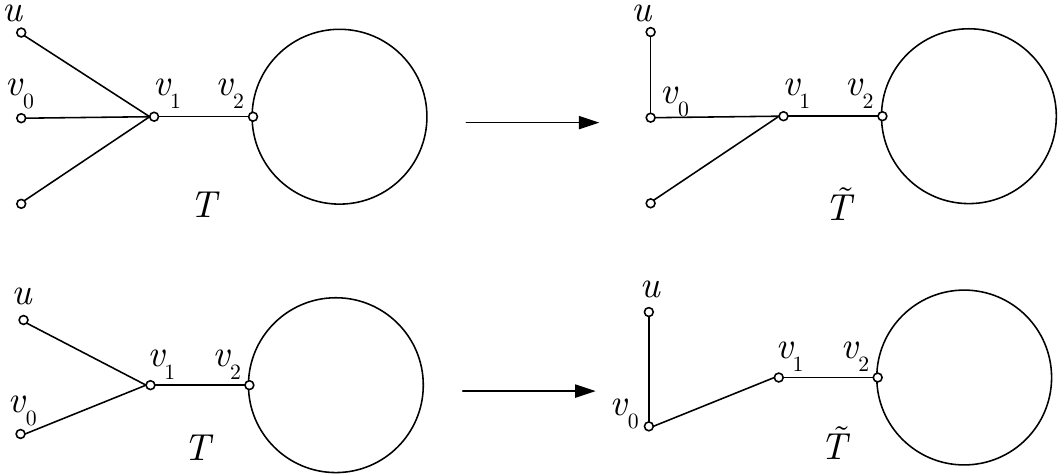}\\
  \caption{ Structure of $T$ in Claim 1.}
  \label{case1}
\end{figure}

On the contrary, suppose $3\le {\rm deg}_T(v_{1})\le 4$. By the choice of $P$, all components of $T-v_1$ are isolated vertices except the component that contains $v_2$. Let $u$ be a  pendant vertex that is adjacent to $v_1$ where $u\not=v_0$.
Let ${\widetilde{T}  =T+v_0u-v_1u}$, see Figure \ref{case1}. By Lemma \ref{lemma2.6}, $d(T)< d(\widetilde{T} )$. Furthermore, $\widetilde{T} \ncong F_{n} $ (by Theorem \ref{theorem2.4}). It follows that $d(T)< d(\widetilde{T})<f(n)$. This is contradictory because $T$ is the tree with the second-largest number of dissociation sets. Thus, $ {\rm deg}_T(v_{1})=2$.

Let $C$ be an arbitrary component of $T-v_2$ such that $V(P)\cap V(C)=\emptyset$. Then there are only two possible structures for $C$. Specifically, we have the following claim.

{\bf Claim 2: } If  $C$ is not an isolated vertex, then $C\cong K_{1,1} $.

In fact, suppose that $C$ is not an isolated vertex. By the choice of $P$, the diameter of the induced subgraph $T\left [ V(C)\cup \left \{ v_{2} \right \}  \right ] $ is no more than two. Thus, $C\cong K_{1,t} $. By Corollary \ref{corollary2.7}, $t\le3$.
To complete the proof of the claim, we need to show $t=1$. In fact, suppose $t\ge 2$. Choose a path $uvw$ from $C$. Clearly, $v$ is the center of $C$ which is adjacent to $v_2$ in $T$.  By Lemma \ref{lemma2.6} and Theorem \ref{theorem2.4}, $d(T)<d(T+wu-uv)<f(n)$. This contradicts the fact that $T$ is the tree with the second-largest number of dissociation sets. It follows that $C\cong K_{1,1}$.

\begin{figure}[htbp!]
  \centering
  \includegraphics[width=0.4\textwidth]{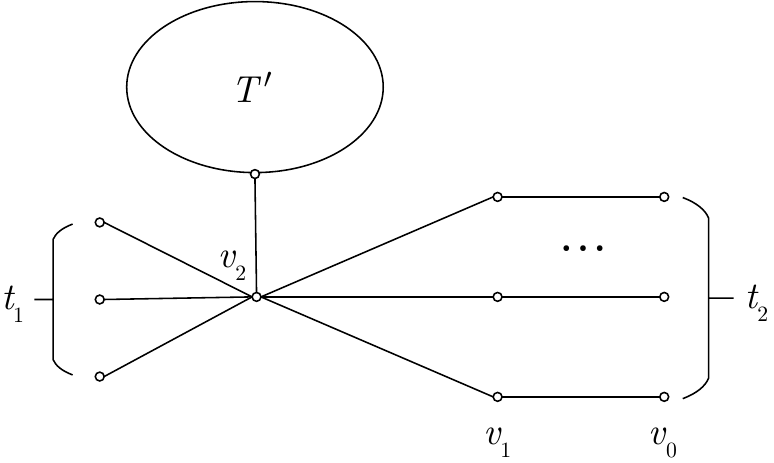}\\
  \caption{Structure of $T$ in Claim 2.}
  \label{case3(1)}
\end{figure}

Thus, every component of $T-v_2$ that does not contain $v_3$ is either an isolated vertex or an isolated edge. Suppose that there are $t_1$ isolated vertices and $t_2$ isolated edges in $T-v_2$, see Figure \ref{case3(1)}.
By Corollary \ref{corollary2.7}, $t_1\le 3$. Let $T'$ be the component of $T-v_2$ that contains $v_3$ and let $|V(T')|=m$.
Since $|V(T)|\ge 11$, $|V(T)|=t_1+2t_2+1+m\ge 11$. It follows that $2t_2+m\ge 7$. Thus,  $m\ge 7$ or $t_2\ge 1$. Because $P=v_0v_1\cdots v_l$ is the longest path in $T$, $l\ge 4$. Specifically, we have the following claim.

{\bf Claim 3:  $l=4$. }

In contrast, suppose that $l>4$.  Then $3\le m\le n-3$. By Lemma \ref{lemma2.2},
\[\begin{aligned}
d(T) & = d(T-v_{0})+d(T-v_{0}-v_{1} )+d(T-v_{0}-v_{1}-v_{2} ).\\
\end{aligned}\]
By Theorem \ref{theorem2.4},
\[\begin{aligned}d(T-v_{0} -v_{1} -v_{2}) & \leq 2^{n-3-m}\cdot f(m) \leq 2^{n-6}\cdot f(3) =7\cdot 2^{n-6}  \\\end{aligned}.\]
\added[id=Reviewer 1]{We have found that there are two different cases of $d(T-v_{0})$. Let us discuss them separately below.}

\added[id=Reviewer 1]{\bf Case 1:  $d(T-v_0)=f(n-1)$.}

\begin{figure}[htbp!]
  \centering
  \includegraphics[width=0.6\textwidth]{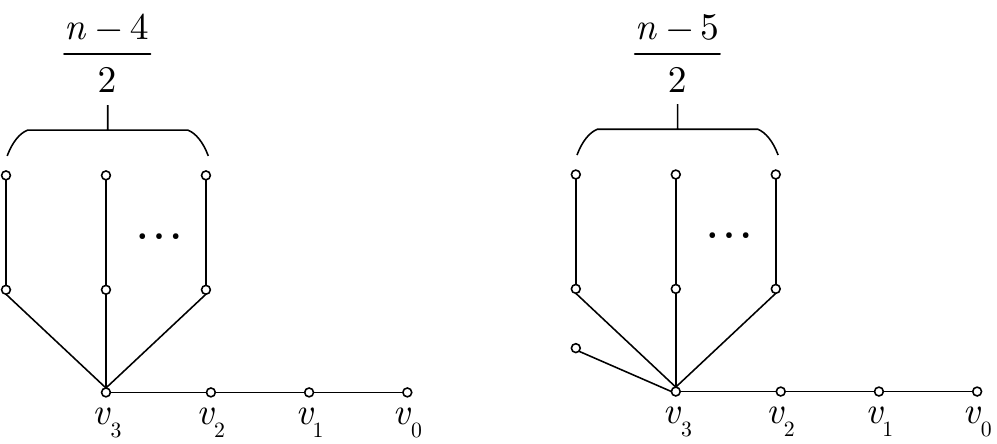}\\
  \caption{ The two possible structures of $T$ if $d(T-v_0)=f(n-1)$.}
  \label{subcase}
\end{figure}

\hspace{2.0em}  By Theorem \ref{theorem2.4}, $T-v_0\cong F_{n-1}$. It follows that $t_1=0$, $t_2=1$.
Thus, there are only two possible structures for  $T$ as shown in Figure \ref{subcase}.
Suppose that $n$ is even \added[id=Reviewer 2]{($n\ge 11$)}.
By Lemma \ref{lemma2.1},  \[\begin{aligned}
d(T) & = d(T-v_3)+d(T-N[v_3])+\sum_{u\in N(v_3)}d(T-N[v_3]\cup N[u]) \\
&= 7\cdot 2^{n-4} +4\cdot 2^{\frac{n-4}{2} } +(\frac{n-4}{2} \cdot 4\cdot 2^{\frac{n-4}{2}-1} + 2\cdot 2^{\frac{n-4}{2} }) \\
&=7\cdot 2^{n-4}+(n+2)\cdot 2^{\frac{n-4}{2} }\\
&< h(n)=2^{n-1}  +(n+12)\cdot 2^{\frac{n-8}{2}}.
\end{aligned}\]
This is contradictory.

\hspace{2.0em}  Suppose $n$ is odd ($n\ge 11$). By Lemma \ref{lemma2.1},
 \[\begin{aligned}
d(T) & = d(T-v_3)+d(T-N[v_3])+\sum_{u\in N(v_3)}d(T-N[v_3]\cup N[u]) \\
&= 2\cdot 7\cdot 2^{n-5}  +4\cdot 2^{\frac{n-5}{2} }+(4\cdot 2^{\frac{n-5}{2} }+\frac{n-5}{2} \cdot 4\cdot 2^{\frac{n-5}{2}-1} + 2\cdot 2^{\frac{n-5}{2} }) \\
&=7\cdot 2^{n-4}+(n+5)\cdot 2^{\frac{n-5}{2} }\\
&< h(n)=2^{n-1}  +(n+9)\cdot 2^{\frac{n-7}{2} }.
\end{aligned}\]
This is contradictory.

\added[id=Reviewer 1]{\bf Case 2:  $d(T-v_0)< f(n-1)$.}

\hspace{2.0em} Note that $T-v_0$ is connected and $d(T-v_0)< f(n-1)$. By induction hypothesis, $d(T-v_0)\le h(n-1)$. Note that $T-v_0-v_1$ is connected. By Theorem \ref{theorem2.4}, $d(T-v_0-v_1)\le f(n-2)$.

\hspace{2.0em} Suppose that $n$ is even ($n\ge 11$). By Lemma \ref{lemma2.1},  \[\begin{aligned}
d(T) &= d(T-v_0)+d(T-v_0-v_1)+d(T-v_0-v_1-v_2) \\
 &\le h(n-1)+f(n-2)+7\cdot 2^{n-6} \\
& = [2^{n-2}+(n+8)\cdot 2^{\frac{n-8}{2} }]+[2^{n-3}+(n+4)\cdot 2^{\frac{n-8}{2} }] +7\cdot 2^{n-6}  \\
&=2^{n-1}  +(n+12)\cdot 2^{\frac{n-8}{2} }-(2^{n-6}-n\cdot 2^{\frac{n-8}{2} })\\
&< h(n)=2^{n-1}  +(n+12)\cdot 2^{\frac{n-8}{2} }.
\end{aligned}\]
The last inequality follows from the fact that $2^{n-6} -n\cdot2^{\frac{n-8}{2} }>0$ when $n\ge 12 $.  This is contradictory.

\hspace{2.0em} Suppose that $n$ is odd ($n\ge 11$). By Lemma \ref{lemma2.1}, \[\begin{aligned}
d(T) &= d(T-v_0)+d(T-v_0-v_1)+d(T-v_0-v_1-v_2) \\
&\le h(n-1)+f(n-2)+7\cdot 2^{n-6} \\
& = [2^{n-2}+(n+11)\cdot 2^{\frac{n-9}{2} }]+[2^{n-3}+(n+1)\cdot 2^{\frac{n-7}{2} } ]+7\cdot 2^{n-6}  \\
&=2^{n-1}  +(n+9)\cdot 2^{\frac{n-7}{2} }-(2^{n-6}-(n-5)\cdot 2^{\frac{n-9}{2} })\\
&< h(n)=2^{n-1}  +(n+9)\cdot 2^{\frac{n-7}{2} }.\end{aligned}\]
The last inequality follows from the fact that $2^{n-6} -(n-5)\cdot 2^{\frac{n-9}{2} }>0$ when $n\ge 11 $.  This is contradictory. Here, we have completed the proof of Claim $3$.

\begin{figure}[htbp!]
  \centering
  \includegraphics[width=0.8\textwidth]{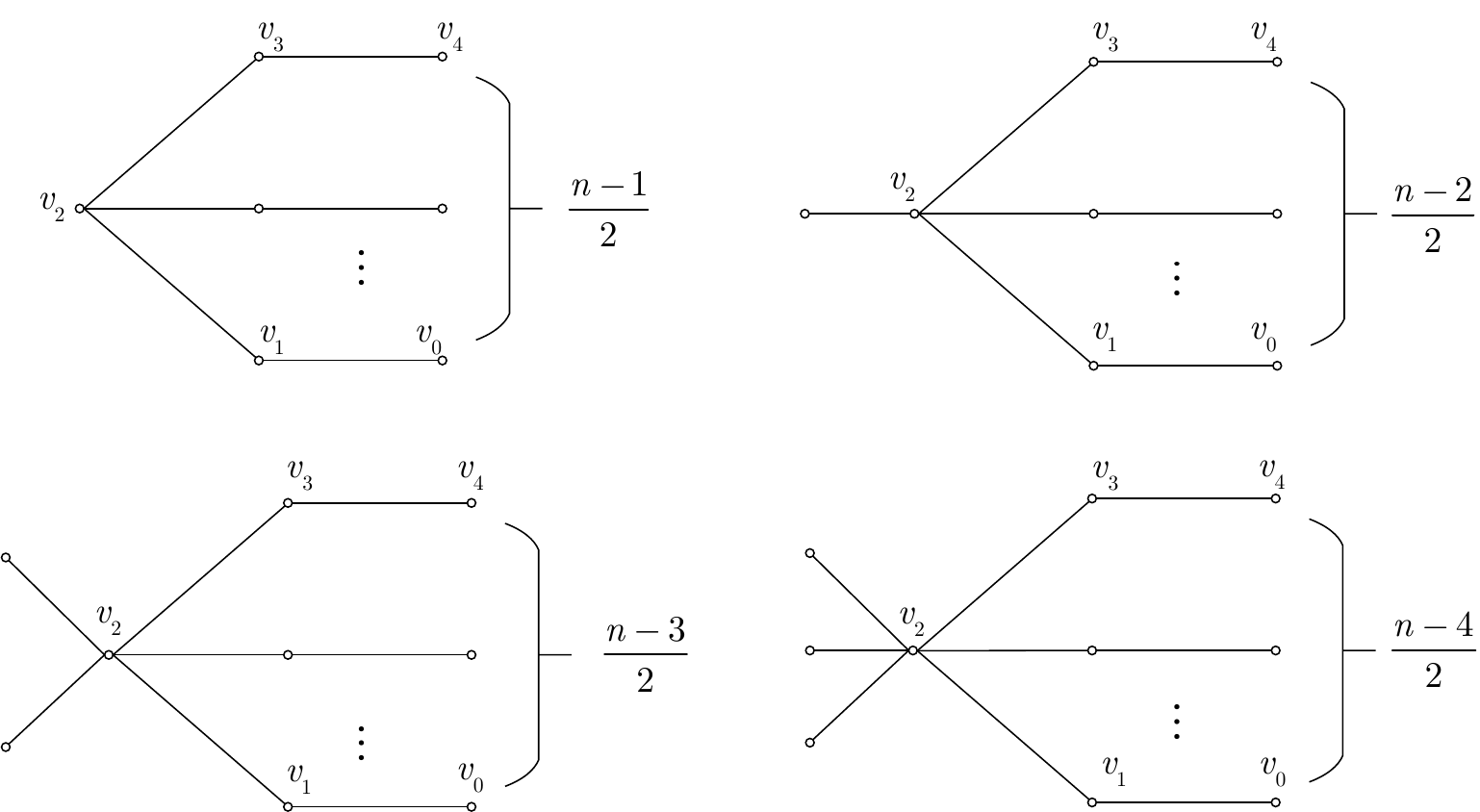}\\
  \caption{ The possible structure of $T$. }
  \label{case3(2)}
\end{figure}

By Claim $3$, the component of $T-v_2$ that contains $v_3$ is a star $K_{1, x}$.  Relabel the vertices of the path $P=v_0v_1\cdots v_l$ as $u_0u_1\cdots u_l$ where $u_i=v_{l-i}$ for all $i=0, 1, \cdots, l$. By Claim $1$, $x=1$. Thus, every component of $T-v_2$ is $K_1$ or $K_2$.  It follows that  
\begin{center}$T\in \{K_{1} \ast \frac{n-1}{2} K_{2}, K_{2} \ast \frac{n-2}{2} K_{2}, K_{2} \ast (\frac{n-3}{2} K_{2} \cup K_{1}), K_{2} \ast (\frac{n-4}{2} K_{2} \cup 2K_{1} )\}$(see Figure \ref{case3(2)}),  \end{center}
the integer $n$ in the first and third elements of the set above are odd, and the other two are even.  
By Theorem \ref{theorem2.4}, $K_{1} \ast \frac{n-1}{2} K_{2}$ and $ K_{2} \ast \frac{n-2}{2}K_{2}$ are two trees with the largest number of dissociation sets. Hence,  $T\in \{ K_{2} \ast (\frac{n-3}{2} K_{2} \cup K_{1} ), K_{2} \ast (\frac{n-4}{2} K_{2} \cup 2K_{1} )\}$, The integer $n$ in the first element of the set is odd, and the second is even.
\hfill$\Box$

\begin{theorem}
Let $G$ be a connected graph of order $n\ge 9$. If $G\ncong F_{n} $, then $d(G)\le h(n)$, the equality holds if and only if $G\in\left \{U_{n}, T_{n} \right \} $, where $F_{n} $, $U_{n}$, $h(n)$, and $T_{n}$ are defined in Theorem \ref{theorem2.4}, Theorem \ref{theorem2.5}, and Theorem \ref{theorem3.1} respectively.
\end{theorem}

\noindent{\bf Proof.} Let $c(G)$ be the dimension of cycle spaces of $G$. If $c(G)\le 1$,  then $G$ is a tree or a unicyclic graph. By Theorem \ref{theorem2.5} and Theorem \ref{theorem3.1}, the result holds. To complete the proof, we just need to show that $d(G)<h(n)$ if $c(G)=d\ge 2$.

In fact, suppose $d\ge 2$. Since $G$ is connected, we can choose $d-1$ edges $u_1v_1, u_2v_2, \ldots, u_{d-1} v_{d-1} $ of $G$ such that the subgraph
\[
G - u_1v_1- u_2v_2- \cdots- u_{d-1} v_{d-1} := {G}'
\]
is a unicyclic graph. Then, we get \[\begin{aligned}
d(G) & \leq d(G-u_{1} v_{1})\\&\leq \cdots \\
& \leq d(G-u_{1} v_{1}-u_{2} v_{2}-\cdots-u_{d-2} v_{d-2}) \\
& \leq d(G-u_{1} v_{1}-u_{2} v_{2}-\cdots-u_{d-2} v_{d-2}-u_{d-1} v_{d-1})\\&=d({G}') \leq h(n).
\end{aligned}\]
Hence, $d(G)=h(n)$ if and only if all of the above inequalities hold.  By Theorem \ref{theorem2.5}, $d({G}')=h(n)$ \added[id=Reviewer 2]{if and only if} ${G}' \cong U _{n} $. By Lemma \ref{lemma2.3}, $d(G-u_{1} v_{1}-u_{2} v_{2}-\cdots-u_{d-2} v_{d-2})=d({G}')$ is equivalent to \added[id=Reviewer 2]{$u_{d-1}$ and $v_{d-1}$} are true twins in $G-u_{1} v_{1}-u_{2} v_{2}-\cdots-u_{d-2} v_{d-2}$\added[id=Reviewer 2]{.} \added[id=Reviewer 2]{It follows that} \added[id=Reviewer 2]{$u_{d-1}$ and $v_{d-1}$} are false twins in ${G}'$. But this contradicts the fact that there are no false twins in ${G}' \cong U _{n} $, where $U _{n}$ is defined in Theorem \ref{theorem2.5}. Then $d(G)<h(n)$. 
\hfill$\Box$

\bibliographystyle{plain}
\bibliography{main}

\newpage
\section*{Appendix}
\label{Appendix}

\begin{figure}[htbp!]
  \centering
  \includegraphics[width=0.8\textwidth]{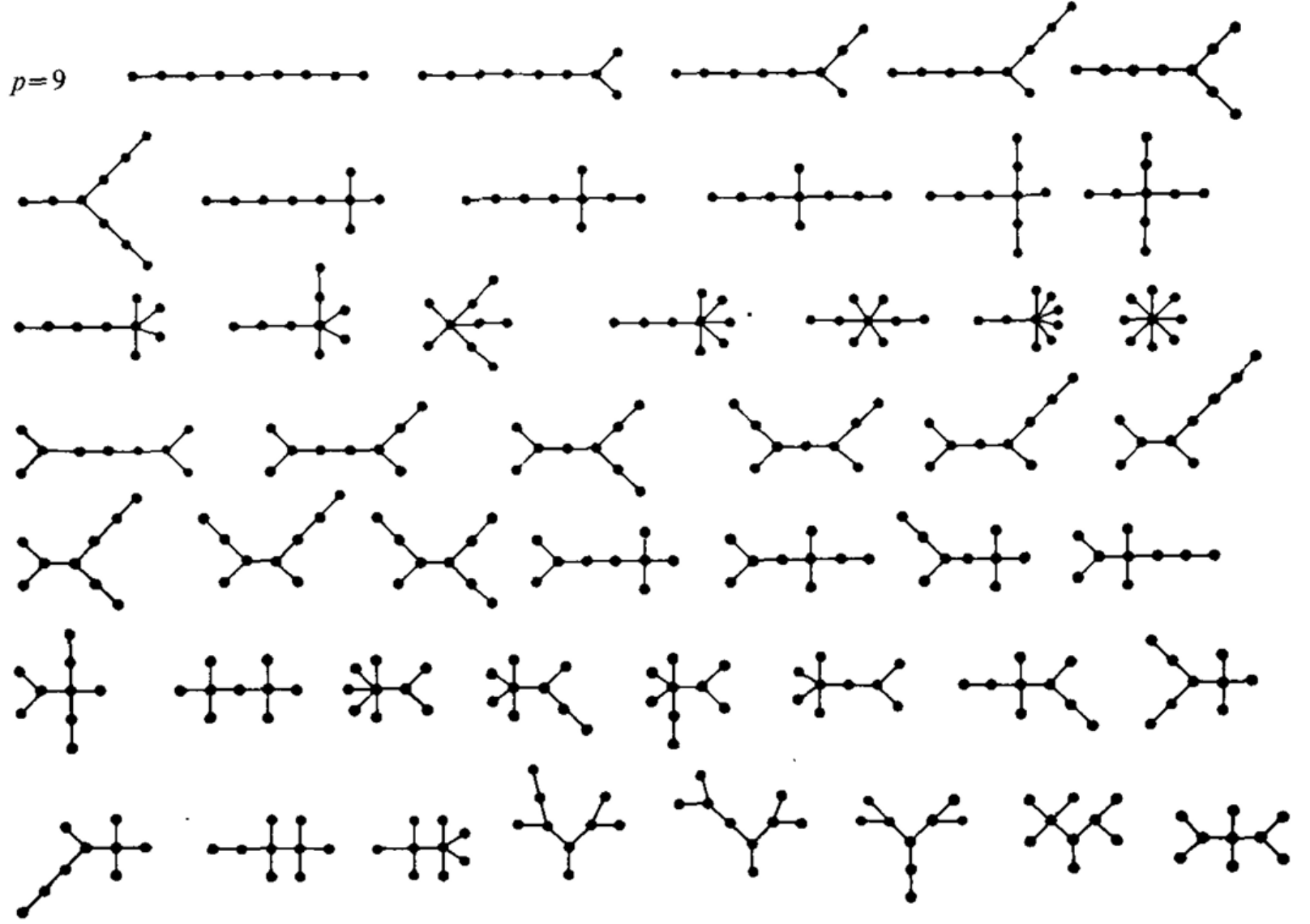}\\
  \caption{ When $n=9$, 47 structures of tree $T$. }
  \label{n=9}
\end{figure}

In the appendix, we will prove Lemma \ref{lemma3.1}.  When $n=9$, there are $47$ different trees in the sense of isomorphism, we label these trees as $T_{1}^{9} ,T_{2}^{9},\cdots T_{47}^{9}$ in order, \footnote{This figure is from page 233 of \cite{Harary1969}.} see Figure~\ref{n=9}.

By Lemma \ref{lemma2.6}, it is easy to see that

\begin{center}
$\begin{array}{llll}
d(T_{2}^{9} )<d(T_{1}^{9}),& d(T_{7}^{9} )<d(T_{3}^{9}),& d(T_{8}^{9} )<d(T_{5}^{9}),& d(T_{9}^{9} )<d(T_{6}^{9}),\vspace{1.0ex}\\
d(T_{13}^{9} )<d(T_{10}^{9}),& d(T_{19}^{9} )<d(T_{1}^{9}),& d(T_{20}^{9} )<d(T_{3}^{9}),& d(T_{21}^{9} )<d(T_{5}^{9}),\vspace{1.0ex}\\
d(T_{23}^{9} )<d(T_{4}^{9}),& d(T_{24}^{9} )<d(T_{4}^{9}),& d(T_{25}^{9} )<d(T_{6}^{9}),& d(T_{28}^{9} )<d(T_{3}^{9}),\vspace{1.0ex}\\
d(T_{29}^{9} )<d(T_{5}^{9}),& d(T_{30}^{9} )<d(T_{22}^{9}),& d(T_{31}^{9} )<d(T_{6}^{9}),& d(T_{32}^{9} )<d(T_{10}^{9}),\vspace{1.0ex}\\
d(T_{33}^{9} )<d(T_{22}^{9}),& d(T_{36}^{9} )<d(T_{10}^{9}),& d(T_{38}^{9} )<d(T_{27}^{9}),& d(T_{39}^{9} )<d(T_{27}^{9}),\vspace{1.0ex}\\
d(T_{40}^{9} )<d(T_{26}^{9}),& d(T_{41}^{9} )<d(T_{27}^{9}),& d(T_{43}^{9} )<d(T_{26}^{9}),& d(T_{44}^{9} )<d(T_{4}^{9}),\vspace{1.0ex}\\
d(T_{45}^{9} )<d(T_{6}^{9}),& d(T_{46}^{9} )<d(T_{26}^{9}),& d(T_{47}^{9} )<d(T_{6}^{9}).
\end{array}$
\end{center}
By Corollary \ref{corollary2.7}, $T_{12}^{9}$, $T_{15}^{9}$, $T_{16}^{9}$, $T_{17}^{9}$, $T_{18}^{9}$, $T_{34}^{9}$, $T_{35}^{9}$, $T_{37}^{9}$, $T_{42}^{9}$ are not the trees with the second-largest number of dissociation sets. By Lemma \ref{lemma2.2} , Theorem \ref{theorem2.4} and Theorem \ref{theorem2.5}, we have $d(T_{1}^{9} )=274<h(9)=292$, $d(T_{11}^{9} )=f(9)=304$, $d(T_{14}^{9} )=h(9)=292$. Therefore, the trees with the second-largest number of dissociation sets belong to $\{T_{3}^{9},T_{4}^{9},T_{5}^{9},T_{6}^{9},T_{10}^{9},T_{22}^{9},T_{26}^{9},T_{27}^{9}\}$.

\begin{figure}[htbp!]
  \centering
  \includegraphics[width=0.8\textwidth]{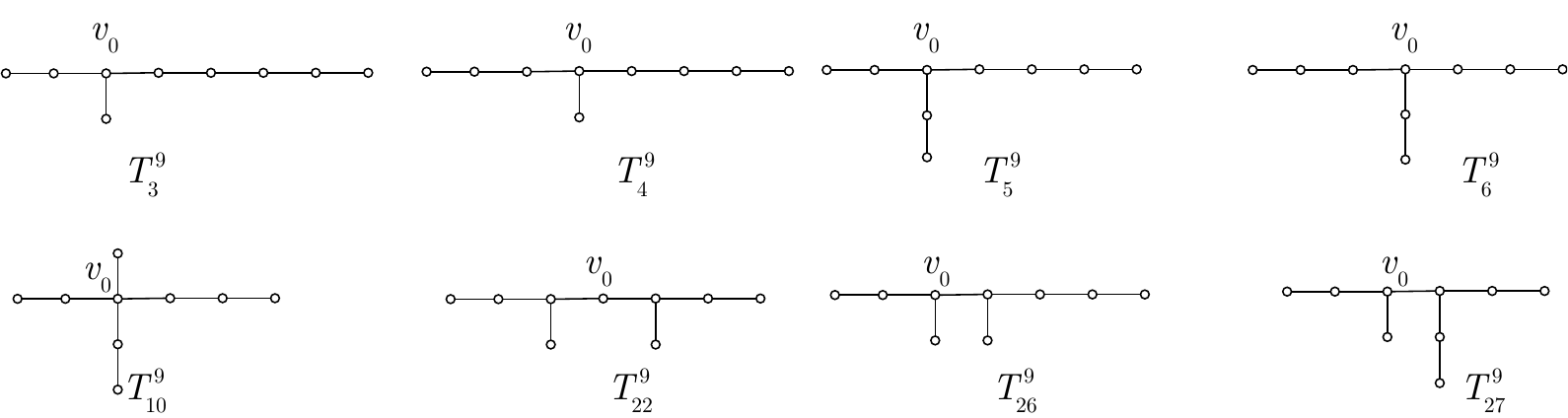}\\
  \caption{ Eight trees that may have the second-largest number of dissociation sets. }
  \label{n=9(1)}
\end{figure}

By Lemma \ref{lemma2.2},  $d(P_n)=d(P_{n-1})+d(P_{n-2})+d(P_{n-3})$. It can be easily calculated that $d(P_1)=2$, $d(P_2)=4$, $d(P_3)=7$, $d(P_4)=13$, $d(P_5)=24$, $d(P_6)=44$, $d(P_7)=81$, $d(P_8)=149$, $d(P_9)=274$. By Lemma \ref{lemma2.2} , we can calculate the number of dissociation sets of these eight graphs as follows by choosing $v_0$ as shown in Figure \ref{n=9(1)}.

In fact,
\[\begin{aligned}
d(T_{3}^{9} ) &= d(T_{3}^{9}-v_0)+d(T_{3}^{9}-N[v_0])+\sum_{u\in N(v_0)} d(T_{3}^{9}-N[v_0]\cup N[u]) \\
& = d(P_2)\cdot d(P_1)\cdot d(P_5)  +2\cdot d(P_1)\cdot d(P_4)+d(P_4)+d(P_1)\cdot d(P_3)\\
&=271<h(9);\\
\end{aligned}\]
\[\begin{aligned}
d(T_{4}^{9} ) &= d(T_{4}^{9}-v_0)+d(T_{4}^{9}-N[v_0])+\sum_{u\in N(v_0)} d(T_{4}^{9}-N[v_0]\cup N[u]) \\
& = d(P_3)\cdot d(P_4)\cdot d(P_1)  +2\cdot d(P_2)\cdot d(P_3)+d(P_1)\cdot d(P_3)+d(P_2)^{2} \\
&=268<h(9);\\
d(T_{5}^{9} ) &= d(T_{5}^{9}-v_0)+d(T_{5}^{9}-N[v_0])+\sum_{u\in N(v_0)} d(T_{5}^{9}-N[v_0]\cup N[u]) \\
& = d(P_2)^{2}\cdot  d(P_4)  +d(P_1)^{2} \cdot d(P_3)+d(P_1)\cdot d(P_3)\cdot 2+d(P_1)^{2}\cdot d(P_2)\\
&=280<h(9);\\
d(T_{6}^{9} ) &= d(T_{6}^{9}-v_0)+d(T_{6}^{9}-N[v_0])+\sum_{u\in N(v_0)} d(T_{6}^{9}-N[v_0]\cup N[u]) \\
& = d(P_2)\cdot d(P_3)^{2}  +d(P_1)\cdot d(P_2)^{2}+d(P_2)\cdot d(P_1)^{2}\cdot 2+d(P_2)^{2}\\
&=276<h(9);\\
d(T_{10}^{9} ) &= d(T_{10}^{9}-v_0)+d(T_{10}^{9}-N[v_0])+\sum_{u\in N(v_0)} d(T_{10}^{9}-N[v_0]\cup N[u]) \\
& = d(P_1)\cdot d(P_2)^{2}\cdot d(P_3)  +d(P_2)\cdot d(P_1)^{2}\cdot 2+2\cdot d(P_1)\cdot d(P_2)+d(P_1)^{3}\\
&=280<h(9);\\
\end{aligned}\]
\[\begin{aligned}
d(T_{22}^{9} ) &= d(T_{22}^{9}-v_0)+d(T_{22}^{9}-N[v_0])+\sum_{u\in N(v_0)} d(T_{22}^{9}-N[v_0]\cup N[u]) \\
& = d(P_4)^{2}+d(P_1)^{2}\cdot d(P_2)^{2}+2\cdot d(P_1)^{2}\cdot d(P_2)\\
&=265<h(9);\\
d(T_{26}^{9} ) &= d(T_{26}^{9}-v_0)+d(T_{26}^{9}-N[v_0])+\sum_{u\in N(v_0)} d(T_{26}^{9}-N[v_0]\cup N[u]) \\
& = d(P_2)\cdot d(P_1)\cdot d(P_5)  +2\cdot d(P_1)^{2}\cdot d(P_3)+d(P_1)\cdot d(P_3)+d(P_1)\cdot d(P_2)\\
&=270<h(9);\\
d(T_{27}^{9} ) &= d(T_{27}^{9}-v_0)+d(T_{27}^{9}-N[v_0])+\sum_{u\in N(v_0)} d(T_{27}^{9}-N[v_0]\cup N[u]) \\
& = d(P_2)\cdot d(P_1)\cdot d(P_5)  +2\cdot d(P_1)\cdot d(P_2)^{2}+d(P_2)^{2}+d(P_1)^{3}\\
&=280<h(9).\\
\end{aligned}\]
So, when $n=9$, the conclusion is true.

\begin{figure}[htbp!]
  \centering
  \includegraphics[width=0.9\textwidth]{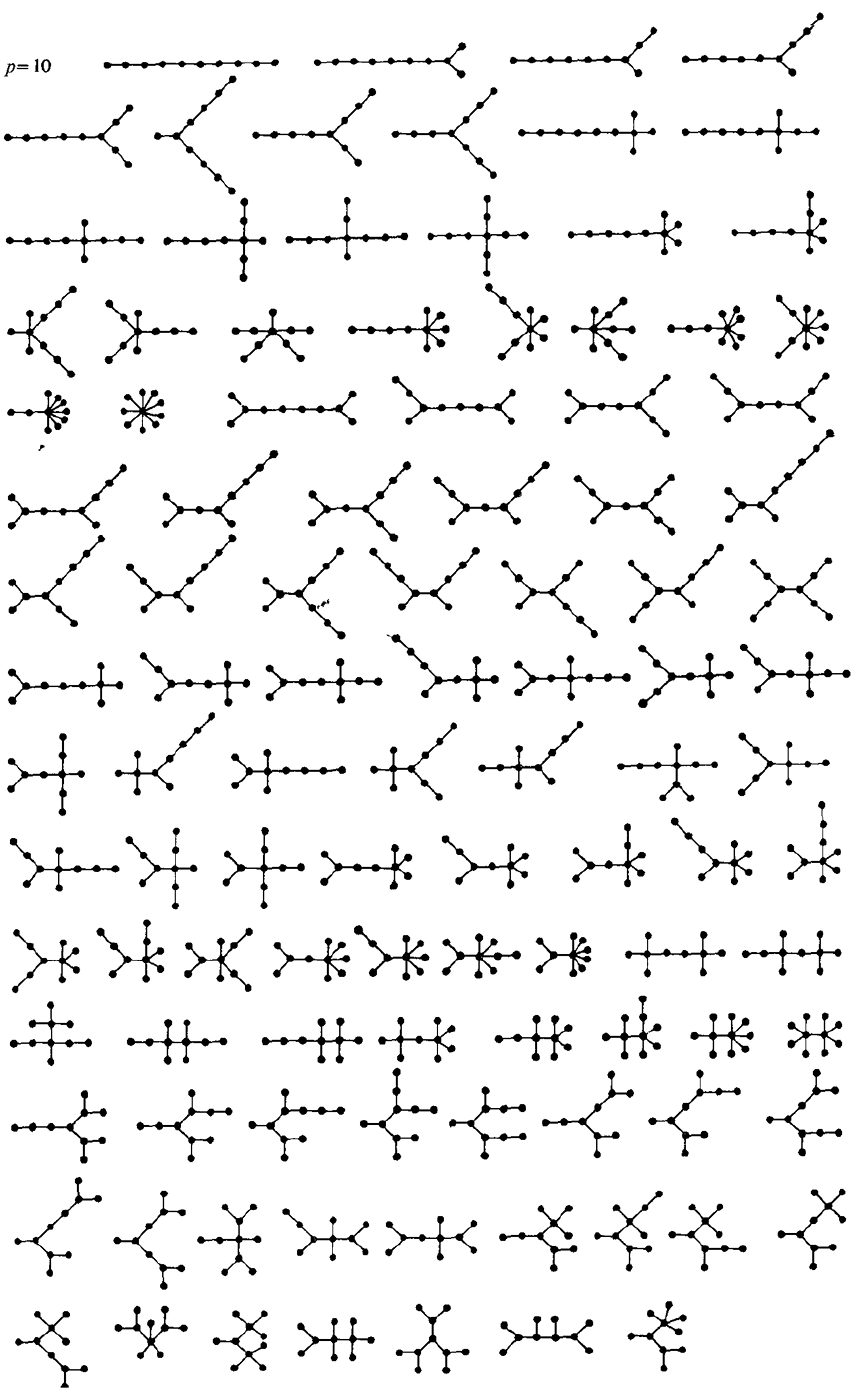}\\
  \caption{ When $n=10$, 106 structures of tree $T$. }
  \label{n=10}
\end{figure}

When $n=10$, there are $106$ different trees in the sense of isomorphism, we label these trees as $T_{1}^{10} ,T_{2}^{10}$, $\cdots$, $T_{106}^{10}$ in order, \footnote{This figure is from page 234 of \cite{Harary1969}.} see Figure \ref{n=10}. By Lemma \ref{lemma2.6}, it is easy to see that 

\begin{figure}[htbp!]
  \centering
  \includegraphics[width=0.9\textwidth]{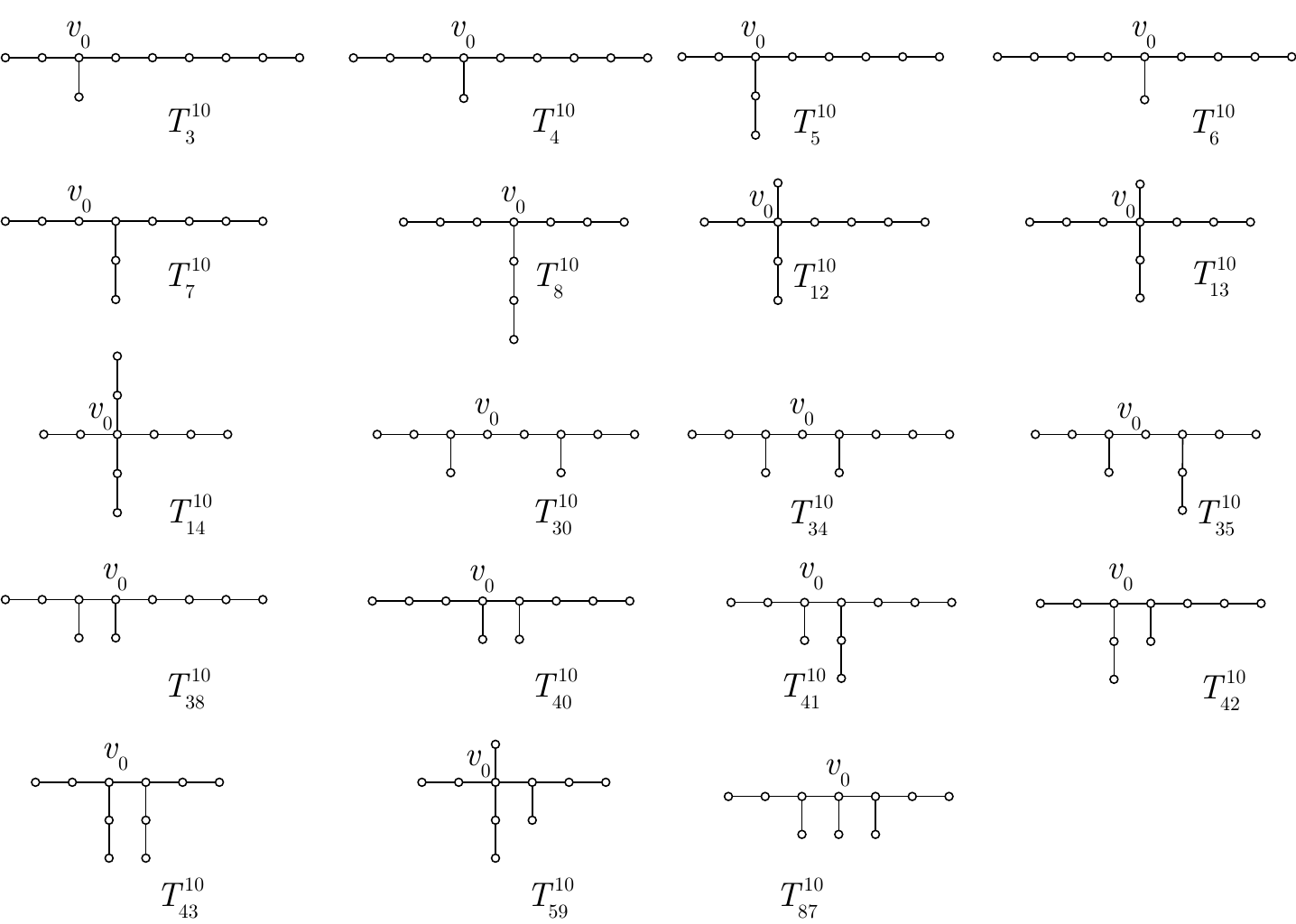}\\
  \caption{ Nineteen trees that may have the second-largest number of dissociation sets. }
  \label{n=10(1)}
\end{figure}

\begin{center}
$\begin{array}{llll}
  d(T_{2}^{10} )<d(T_{1}^{10}), & d(T_{9}^{10} )<d(T_{3}^{10}), & d(T_{10}^{10} )<d(T_{5}^{10}), &d(T_{11}^{10} )<d(T_{7}^{10}),\vspace{1.0ex} \\
  d(T_{16}^{10} )<d(T_{12}^{10}), & d(T_{17}^{10} )<d(T_{13}^{10}),&d(T_{18}^{10} )<d(T_{14}^{10}),&d(T_{27}^{10} )<d(T_{1}^{10}),\vspace{1.0ex}\\
  d(T_{28}^{10} )<d(T_{3}^{10}), & d(T_{29}^{10} )<d(T_{5}^{10}),&d(T_{31}^{10} )<d(T_{4}^{10}),&d(T_{32}^{10} )<d(T_{6}^{10}),\vspace{1.0ex}\\
  d(T_{33}^{10} )<d(T_{7}^{10}), & d(T_{36}^{10} )<d(T_{4}^{10}),&d(T_{37}^{10} )<d(T_{7}^{10}),&d(T_{39}^{10} )<d(T_{8}^{10}),\vspace{1.0ex}\\
  d(T_{44}^{10} )<d(T_{3}^{10}), & d(T_{45}^{10} )<d(T_{30}^{10}),&d(T_{46}^{10} )<d(T_{5}^{10}),&d(T_{47}^{10} )<d(T_{34}^{10}),\vspace{1.0ex}\\
  d(T_{48}^{10} )<d(T_{7}^{10}), & d(T_{49}^{10} )<d(T_{35}^{10}),&d(T_{50}^{10} )<d(T_{35}^{10}),&d(T_{51}^{10} )<d(T_{12}^{10}),\vspace{1.0ex}\\
  d(T_{52}^{10} )<d(T_{38}^{10}), & d(T_{53}^{10} )<d(T_{7}^{10}),&d(T_{54}^{10} )<d(T_{41}^{10}),&d(T_{55}^{10} )<d(T_{42}^{10}),\vspace{1.0ex}\\
  d(T_{56}^{10} )<d(T_{13}^{10}), & d(T_{57}^{10} )<d(T_{43}^{10}),&d(T_{58}^{10} )<d(T_{41}^{10}),&d(T_{60}^{10} )<d(T_{14}^{10}),\vspace{1.0ex}\\
  d(T_{63}^{10} )<d(T_{12}^{10}), & d(T_{65}^{10} )<d(T_{13}^{10}),&d(T_{67}^{10} )<d(T_{59}^{10}),&d(T_{68}^{10} )<d(T_{14}^{10}),\vspace{1.0ex}\\
  d(T_{73}^{10} )<d(T_{30}^{10}), & d(T_{74}^{10} )<d(T_{35}^{10}),&d(T_{75}^{10} )<d(T_{59}^{10}),&d(T_{76}^{10} )<d(T_{43}^{10}),\vspace{1.0ex}\\
  d(T_{77}^{10} )<d(T_{41}^{10}), & d(T_{80}^{10} )<d(T_{59}^{10}),&d(T_{83}^{10} )<d(T_{8}^{10}),&d(T_{84}^{10} )<d(T_{41}^{10}),\vspace{1.0ex}\\
  d(T_{85}^{10} )<d(T_{40}^{10}), & d(T_{86}^{10} )<d(T_{42}^{10}),&d(T_{88}^{10} )<d(T_{7}^{10}),&d(T_{89}^{10} )<d(T_{34}^{10}),\vspace{1.0ex}\\
  d(T_{90}^{10} )<d(T_{38}^{10}), & d(T_{91}^{10} )<d(T_{4}^{10}),&d(T_{92}^{10} )<d(T_{6}^{10}),&d(T_{93}^{10} )<d(T_{13}^{10}),\vspace{1.0ex}\\
  d(T_{94}^{10} )<d(T_{41}^{10}), & d(T_{95}^{10} )<d(T_{7}^{10}),&d(T_{96}^{10} )<d(T_{41}^{10}),&d(T_{97}^{10} )<d(T_{42}^{10}),\vspace{1.0ex}\\
  d(T_{98}^{10} )<d(T_{87}^{10}), & d(T_{99}^{10} )<d(T_{34}^{10}),&d(T_{100}^{10} )<d(T_{38}^{10}),&d(T_{101}^{10} )<d(T_{13}^{10}),\vspace{1.0ex}\\
  d(T_{102}^{10} )<d(T_{87}^{10}), & d(T_{103}^{10} )<d(T_{41}^{10}),&d(T_{104}^{10} )<d(T_{8}^{10}),&d(T_{105}^{10} )<d(T_{40}^{10}).
\end{array}$
\end{center}
By Corollary \ref{corollary2.7}, $T_{15}^{10}$, $T_{20}^{10}$, $T_{21}^{10}$, $T_{23}^{10}$, $T_{24}^{10}$, $T_{25}^{10}$, $T_{26}^{10}$, $T_{61}^{10}$, $T_{62}^{10}$, $T_{64}^{10}$, $T_{66}^{10}$, $T_{69}^{10}$, $T_{70}^{10}$, $T_{71}^{10}$, $T_{72}^{10}$, $T_{78}^{10}$, $T_{79}^{10}$, $T_{81}^{10}$, $T_{82}^{10}$, $T_{106}^{10}$ are not the trees with the second-largest number of dissociation sets. By Lemma \ref{lemma2.2} , Theorem \ref{theorem2.4} and Theorem \ref{theorem2.5}, we have $d(T_{1}^{10} )=504<h(10)=556$, $d(T_{19}^{10} )=f(10)=576$, $d(T_{22}^{10} )=h(10)=556$. Therefore, the graphs with the second-largest number of dissociation sets belong to
\begin{center}
$\{T_{3}^{10}$, $T_{4}^{10}$, $T_{5}^{10}$, $T_{6}^{10}$, $T_{7}^{10}$, $T_{8}^{10}$, $T_{12}^{10}$, $T_{13}^{10}$, $T_{14}^{10}$, $T_{30}^{10}$, $T_{34}^{10}$, $T_{35}^{10}$, $T_{38}^{10}$, $T_{40}^{10}$, $T_{41}^{10}$, $T_{42}^{10}$, $T_{43}^{10}$, $T_{59}^{10}$, $T_{87}^{10}\}$.
\end{center}
By Lemma \ref{lemma2.2} , we can calculate the number of dissociation sets of these nineteen trees by choosing $v_0$ as shown in Figure \ref{n=10(1)}.
\[\begin{aligned}
d(T_{3}^{10} ) &= d(T_{3}^{10}-v_0)+d(T_{3}^{10}-N[v_0])+\sum_{u\in N(v_0)} d(T_{3}^{10}-N[v_0]\cup N[u]) \\
& = d(P_2)\cdot d(P_1)\cdot d(P_6)  +2\cdot d(P_1)\cdot d(P_5)+d(P_5)+d(P_1)\cdot d(P_4)\\
&=498<h(10);\\
d(T_{4}^{10} ) &= d(T_{4}^{10}-v_0)+d(T_{4}^{10}-N[v_0])+\sum_{u\in N(v_0)} d(T_{4}^{10}-N[v_0]\cup N[u]) \\
& = d(P_3)\cdot d(P_1)\cdot d(P_5)  +2\cdot d(P_2)\cdot d(P_4)+d(P_1)\cdot d(P_4)+d(P_2)\cdot d(P_3)\\
&=494<h(10);\\
d(T_{5}^{10} ) &= d(T_{5}^{10}-v_0)+d(T_{5}^{10}-N[v_0])+\sum_{u\in N(v_0)} d(T_{5}^{10}-N[v_0]\cup N[u]) \\
& = d(P_2)^{2}\cdot d(P_5)  +d(P_1)^{2}\cdot d(P_4)  +2\cdot d(P_1)\cdot d(P_4)+d(P_1)^{2}\cdot d(P_3)\\
&=516<h(10);\\
d(T_{6}^{10} ) &= d(T_{6}^{10}-v_0)+d(T_{6}^{10}-N[v_0])+\sum_{u\in N(v_0)} d(T_{6}^{10}-N[v_0]\cup N[u]) \\
& = d(P_4)^{2}\cdot d(P_1)  +d(P_3)^{2}\cdot 2  +2\cdot d(P_2)\cdot d(P_3)\\
&=492<h(10);\\
\end{aligned}\]
\[\begin{aligned}
d(T_{7}^{10} ) &= d(T_{7}^{10}-v_0)+d(T_{7}^{10}-N[v_0])+\sum_{u\in N(v_0)} d(T_{7}^{10}-N[v_0]\cup N[u]) \\
& = d(P_2)\cdot d(P_7)+d(P_1)\cdot d(P_2)\cdot d(P_4) +d(P_2)\cdot d(P_4)+d(P_1)^{2}\cdot d(P_3)\\
&=508<h(10);\\
d(T_{8}^{10} ) &= d(T_{8}^{10}-v_0)+d(T_{8}^{10}-N[v_0])+\sum_{u\in N(v_0)} d(T_{8}^{10}-N[v_0]\cup N[u]) \\
& = d(P_3)^{3}+d(P_2)^{3}+3\cdot d(P_1)\cdot d(P_2)^{2}\\
&=503<h(10);\\
d(T_{12}^{10} ) &= d(T_{12}^{10}-v_0)+d(T_{12}^{10}-N[v_0])+\sum_{u\in N(v_0)} d(T_{12}^{10}-N[v_0]\cup N[u]) \\
& = d(P_1)\cdot d(P_2)^{2}\cdot d(P_4)+d(P_1)^{2}\cdot d(P_3)\cdot 2+d(P_1)\cdot d(P_3)\cdot 2+d(P_1)^{2}\cdot d(P_2)\\
&=516<h(10);\\
d(T_{13}^{10} ) &= d(T_{13}^{10}-v_0)+d(T_{13}^{10}-N[v_0])+\sum_{u\in N(v_0)} d(T_{13}^{10}-N[v_0]\cup N[u]) \\
& = d(P_1)\cdot d(P_3)^{2}\cdot d(P_2)+d(P_2)^{2}\cdot d(P_1)\cdot 2+d(P_2)^{2}+d(P_1)^{2}\cdot d(P_2)\cdot 2\\
&=504<h(10);\\
d(T_{14}^{10} ) &= d(T_{14}^{10}-v_0)+d(T_{14}^{10}-N[v_0])+\sum_{u\in N(v_0)} d(T_{14}^{10}-N[v_0]\cup N[u]) \\
& = d(P_2)^{3}\cdot d(P_3)+d(P_1)^{3}\cdot d(P_2)+d(P_1)^{2}\cdot d(P_2)\cdot 3+d(P_1)^{4}\\
&=544<h(10);\\
\end{aligned}
\]

\[\begin{aligned}
d(T_{30}^{10} ) &= d(T_{30}^{10}-v_0)+d(T_{30}^{10}-N[v_0])+\sum_{u\in N(v_0)} d(T_{30}^{10}-N[v_0]\cup N[u]) \\
& = d(P_4)\cdot 23+d(P_1)\cdot d(P_2)\cdot d(P_4)+d(P_1)\cdot d(P_4)+d(P_1)^{2}\cdot d(P_2)^{2}\\
&=493<h(10);\\
d(T_{34}^{10} ) &= d(T_{34}^{10}-v_0)+d(T_{34}^{10}-N[v_0])+\sum_{u\in N(v_0)} d(T_{34}^{10}-N[v_0]\cup N[u]) \\
& = d(P_4)\cdot d(P_5)+d(P_1)^{2}\cdot d(P_2)\cdot d(P_3)+d(P_1)^{2}\cdot d(P_3)+d(P_1)\cdot d(P_2)^{2}\\
&=484<h(10);\\
d(T_{35}^{10} ) &= d(T_{35}^{10}-v_0)+d(T_{35}^{10}-N[v_0])+\sum_{u\in N(v_0)} d(T_{35}^{10}-N[v_0]\cup N[u]) \\
& = d(P_4)\cdot d(P_5)+d(P_1)\cdot d(P_2)^{3}+d(P_1)\cdot d(P_2)^{2}+d(P_1)^{3}\cdot d(P_2)\\
&=504<h(10);\\
d(T_{38}^{10} ) &= d(T_{38}^{10}-v_0)+d(T_{38}^{10}-N[v_0])+\sum_{u\in N(v_0)} d(T_{38}^{10}-N[v_0]\cup N[u]) \\
& = d(P_1)\cdot d(P_4)^{2}+d(P_1)\cdot d(P_2)\cdot d(P_3)\cdot 2+d(P_1)\cdot d(P_3)+d(P_2)^{2}\cdot d(P_1)\\
&=496<h(10);\\
d(T_{40}^{10} ) &= d(T_{40}^{10}-v_0)+d(T_{40}^{10}-N[v_0])+\sum_{u\in N(v_0)} d(T_{40}^{10}-N[v_0]\cup N[u]) \\
& = d(P_1)\cdot d(P_3)\cdot d(P_5)+d(P_1)\cdot d(P_2)\cdot d(P_3)\cdot 2+d(P_1)^{2}\cdot d(P_3)+d(P_2)^{2}\\
&=492<h(10);\\
d(T_{41}^{10} ) &= d(T_{41}^{10}-v_0)+d(T_{41}^{10}-N[v_0])+\sum_{u\in N(v_0)} d(T_{41}^{10}-N[v_0]\cup N[u]) \\
& = d(P_1)\cdot d(P_2)\cdot d(P_6)+d(P_1)\cdot d(P_2)\cdot d(P_3)\cdot 2+d(P_2)\cdot d(P_3)+d(P_1)^{2}\cdot d(P_2)\\
&=508<h(10);\\
d(T_{42}^{10} ) &= d(T_{42}^{10}-v_0)+d(T_{42}^{10}-N[v_0])+\sum_{u\in N(v_0)} d(T_{42}^{10}-N[v_0]\cup N[u]) \\
& = d(P_2)^{2}\cdot d(P_5)+d(P_1)^{3}\cdot d(P_3)+d(P_1)^{2}\cdot d(P_3)\cdot 2+d(P_1)^{2}\cdot d(P_2)\\
&=512<h(10);\\
d(T_{43}^{10} ) &= d(T_{43}^{10}-v_0)+d(T_{43}^{10}-N[v_0])+\sum_{u\in N(v_0)} d(T_{43}^{10}-N[v_0]\cup N[u]) \\
& = d(P_2)^{2}\cdot d(P_5)+d(P_1)^{2}\cdot d(P_2)^{2}+d(P_1)\cdot d(P_2)^{2}\cdot 2+d(P_1)^{4}\\
&=528<h(10);\\
d(T_{59}^{10} ) &= d(T_{59}^{10}-v_0)+d(T_{59}^{10}-N[v_0])+\sum_{u\in N(v_0)} d(T_{59}^{10}-N[v_0]\cup N[u]) \\
& = d(P_1)\cdot d(P_2)^{2}\cdot d(P_4)+d(P_1)^{3}\cdot d(P_2)\cdot 2+d(P_1)^{2}\cdot d(P_2)\cdot 2+d(P_1)^{3}\\
&=520<h(10);\\
d(T_{87}^{10} ) &= d(T_{87}^{10}-v_0)+d(T_{87}^{10}-N[v_0])+\sum_{u\in N(v_0)} d(T_{87}^{10}-N[v_0]\cup N[u]) \\
& = d(P_1)\cdot d(P_4)^{2}+d(P_1)^{2}\cdot d(P_2)^{2}\cdot 2+d(P_1)^{2}\cdot d(P_2)\cdot 2\\
&=498<h(10).\\
\end{aligned}
\]
So the conclusion is true when $n=10$. $\qed$

\end{document}